\documentclass[letterpaper, 10 pt, conference]{ieeeconf}

\IEEEoverridecommandlockouts
\usepackage{cite}
\usepackage{amsmath,amssymb,amsfonts}
\usepackage{amsthm}
\usepackage{algorithmic}
\usepackage{graphicx}
\usepackage{textcomp}
\usepackage{comment}
\usepackage{subcaption}

\usepackage[font=small,labelfont=bf]{caption}

\usepackage{xcolor}
\usepackage[colorlinks = true,
            linkcolor = blue,
            urlcolor  = blue,
            citecolor = blue,
            anchorcolor = blue]{hyperref}

\newtheorem{theorem}{Theorem}

\newtheorem{lemma}{Lemma}
\newtheorem{corollary}{Corollary}

\newtheorem{assumption}{Assumption}

\newtheorem{definition}{Definition}
\newtheorem{example}{Example}

\newtheorem{remark}{Remark}

\usepackage[nameinlink]{cleveref}
\crefname{equation}{}{}
\crefname{theorem}{Theorem}{Theorems}
\crefname{corollary}{Corollary}{Corollaries}
\crefname{example}{Example}{Examples}
\crefname{assumption}{Assumption}{Assumptions}
\crefname{lemma}{Lemma}{Lemmas}
\crefname{proposition}{Proposition}{Propositions}
\crefname{figure}{Figure}{Figures}
\crefname{table}{Table}{Tables}
\crefname{fact}{Fact}{Facts}
\crefname{conjecture}{Conjecture}{Conjectures}
\crefname{section}{Section}{Sections}
\crefname{appendix}{Appendix}{Appendices}
\Crefname{equation}{}{}
\Crefname{theorem}{Theorem}{Theorems}
\Crefname{corollary}{Corollary}{Corollaries}
\Crefname{example}{Example}{Examples}
\Crefname{lemma}{Lemma}{Lemma}
\Crefname{proposition}{Proposition}{Proposition}
\Crefname{figure}{Figure}{Figures}
\Crefname{table}{Table}{Tables}
\Crefname{section}{Section}{Sections}
\Crefname{appendix}{Appendix}{Appendices}

\Crefname{problem}{Problem}{Problem}

\newcommand{\tr}{{{\mathsf T}}}
\newcommand{\Tr}{{{\mathbf{Tr}}}}
\newcommand{\her}{{{\mathsf H}}}

\newcommand{\Hinf}{{\mathcal{H}_\infty}}

\newcommand{\bT}{{\mathbf{T}}}

\newcommand{\diag}{\mathrm{diag}}

\newcommand{\FR}{\mathfrak{R}}
\newcommand{\FI}{\mathfrak{I}}
\newcommand{\FJ}{\mathfrak{J}}

\newcommand{\RR}{\mathbb{R}} 

\begin{document}
\bstctlcite{IEEEexample:BSTcontrol}

\title{\bf 
Policy Optimization in Robust Control: Weak Convexity and Subgradient Methods 
}

\author{Yuto Watanabe, 
Feng-Yi Liao, and Yang Zheng
\thanks{This work is supported by NSF ECCS-2154650, NSF CMMI 2320697, and NSF CAREER 2340713.}
\thanks{Y. Watanabe, F.-Y. Liao,  and Y. Zheng are with the Department of Electrical
and Computer Engineering, University of California San Diego; \texttt{\{y1watanabe,fliao,zhengy\}@ucsd.edu}.}}

\maketitle

\begin{abstract}
Robust control seeks stabilizing policies that perform reliably under adversarial disturbances, with $\mathcal{H}_\infty$ control as a classical formulation. It is known that policy optimization of robust $\mathcal{H}_\infty$ control naturally lead to nonsmooth and nonconvex problems. This paper builds on recent advances in nonsmooth optimization to analyze discrete-time static output-feedback $\mathcal{H}_\infty$ control. We show that the $\mathcal{H}_\infty$ cost is \emph{weakly convex} over any convex subset of a sublevel set. This structural property allows us to establish the first non-asymptotic deterministic convergence rate for the subgradient method under suitable assumptions. In addition, we prove a weak Polyak--\L{}ojasiewicz (PL) inequality in the state-feedback case, implying that all stationary points are globally optimal. We finally present a few numerical examples to validate the theoretical results.
\end{abstract}

\section{Introduction}
\label{sec:introduction}

Ensuring robustness against unknown disturbances is crucial 
in control systems.
In this context, 
$\mathcal{H}_\infty$ control has served as a cornerstone in robust control \cite{zhou1996robust}, which aims to design a stabilizing policy that minimizes the $\mathcal{H}_\infty$ norm of the closed-loop system. 
It has long been recognized that policy optimization for robust control naturally leads to nonconvex and nonsmooth optimization problems \cite{apkarian2006nonsmooth,guo2022global,guo2023complexity,tang2023global,zheng2023benign,zheng2024benign,hu2022connectivity}.
While Riccati equation-based and LMI-based methods have been developed to circumvent these difficulties \cite{doyle1988state,gahinet1994linear}, they are typically limited to unstructured and model-based controller design. 
To address this, some recent works have proposed local search algorithms that directly optimize the original nonconvex nonsmooth problem \cite{apkarian2006nonsmooth,burke2006hifoo,saeki2006fixed,apkarian2009proximity,guo2022global,guo2023complexity}.
In particular,
the authors in \cite{apkarian2006nonsmooth,apkarian2009proximity} have implemented two official MATLAB solvers, \texttt{hinfstruct} and \texttt{systune}, for structured $\mathcal{H}_\infty$ control, with successful real-world applications
(e.g., a space probe by the ESA
\cite{apkarian2017h}). 

However, owing to the inherent nonconvexity and nonsmoothness, the theoretical understanding remains limited, particularly regarding the non-asymptotic convergence of local search algorithms. 
The works \cite{apkarian2006nonsmooth,saeki2006fixed,apkarian2009proximity} proposed \textit{deterministic} local search algorithms with some structure exploitation, but only admit \textit{at most} asymptotic convergence.
As a \textit{probabilistic} approach, \cite{burke2006hifoo}
presented an open-source package \texttt{HIFOO} based on a randomized gradient sampling strategy \cite{burke2005robust}, 
which was also limited to asymptotic guarantees.
Recently, 
\cite{guo2022global,guo2023complexity}
presented algorithms with a non-asymptotic convergence rate in a probabilistic sense
using recent randomization-based techniques in 
\cite{zhang2020complexity,duchi2012randomized}.
However, such probabilistic guarantees are subject to a risk of failure, which may lead to concerns in safety-critical systems.
To our knowledge,
deterministic non-asymptotic convergence rates
are still largely open for $\mathcal{H}_\infty$ optimization.

This work addresses the aforementioned issue of convergence rates via \textit{weak convexity}, and establishes a non-asymptotic convergence rate for the simple subgradient method.
Recent studies \cite{davis2019stochastic,liang2023proximal,liao2024error,liao2025proximal} demonstrate that weakly convex functions 
are arguably one of the broadest classes 
in nonsmooth nonconvex optimization that
admits deterministic non-asymptotic rates,
which cover various practical functions (e.g., convex, smooth, the composite of a convex function and a smooth map).
Building on the recent advances,
we first establish weak convexity of the $\Hinf$ cost
for discrete-time systems with static output-feedback.
While previous works \cite{apkarian2009proximity,noll2013bundle} pointed out a relevant property, called \textit{lower-$C^2$} \cite{rockafellar2009variational},
we demonstrate that the $\Hinf$ cost is more structured than lower-$C^2$.
With this property,
we present a non-asymptotic convergence rate for the subgradient method under the assumption that the cost remains finite.
Further, we present 
a \textit{weak Polyak-\L{}ojasiewicz} (PL) inequality \cite{umenberger2022globally} for the full state-feedback case,
ensuring the global optimality of stationary points despite the nonconvexity.

Our major technical contributions
are threefold:
\begin{enumerate}
    \item \textbf{Weak convexity.} 
Since the stability constraint on feedback gains is inherently nonconvex, the standard definition of weak convexity is not applicable.
Thus, we 
instead
show $m$-weak convexity with some $m>0$ over any convex subset of a sublevel set.
Notably, the constant $m$ is uniform over the sublevel set, which allows for a straightforward extension of the standard weakly convex case \cite{davis2019stochastic};
    \item \textbf{A weak Polyak-Lojasiewicz (PL) inequality.}
    In the full state measurement case, we present a weak PL inequality \cite{umenberger2022globally},
    guaranteeing the global optimality of stationary points. Unlike the differentiability-based result in \cite{umenberger2022globally}, our proof applies to the nonsmooth $\mathcal{H}_\infty$ cost. 
    The argument leverages an analogous idea to the extended convex lifting \cite{zheng2024benign} to analyze the nonsmooth $\mathcal{H}_\infty$ landscape beyond standard stationarity;
    \item \textbf{The subgradient method.}
    Building on weak convexity, we extend the convergence analysis of \cite{davis2019stochastic} to the constrained nonconvex $\mathcal{H}_\infty$ optimization setting. 
    In particular, we show that an $\epsilon$-inexact stationary point can be obtained 
    after $T = \mathcal{O}(1/\epsilon^4)$ iterations, provided the cost remains bounded along the iterates.  To our knowledge, this is the first deterministic complexity guarantee for subgradient methods in robust control.  
\end{enumerate}

The remainder of this paper is organized as follows.
\cref{section:problem-formulation}
formulates the $\Hinf$ optimization problem and provides the problem statement.
Then, \cref{section:landscape} presents
a review of landscape properties and 
our main results,
such as the weak convexity and weak PL condition.
In \cref{section:subgradient_method},
we discuss the subgradient method.
\cref{section:simulations} provides numerical results.
Finally, \cref{section:conclusion} concludes this paper.

\section{Preliminaries and Problem Statement}\label{section:problem-formulation}

\subsection{Robust $\mathcal{H}_\infty$ optimization with static output-feedback}\label{subsection:problem_formulation}

Consider a discrete-time linear time-invariant system\footnote{
We focus on the class of discrete-time systems, since it is 
practical and might also be the simplest case in $\Hinf$ optimization, due to the \textit{coercivity};
see \cite[Lem. 3.2]{guo2023complexity} and \cite[Example 2.3]{zheng2024benign}.
See also
\Cref{remark:continuous-time}.
}
\begin{equation}\label{eq:dynamic}
\begin{aligned}
x_{t+1} & =A x_t+B u_t+B_w w_t,\quad
y_t  =C x_t,
\end{aligned}
\end{equation}
where $x_t \in \mathbb{R}^{n_x}$ is the system state, $u_t \in \mathbb{R}^{n_u}$ is the control input, $y_t \in \mathbb{R}^{n_y}$ is the output measurement, and $w_t \in \mathbb{R}^{n_w}$ is the disturbance, respectively. 
The system matrices are $A \in \mathbb{R}^{n_x \times n_x}, B \in \mathbb{R}^{n_x \times n_u}$, $B_w \in \mathbb{R}^{n_x \times n_w}$, and $ C \in \mathbb{R}^{n_y \times n_x}$. 
The system's initial condition is fixed as $x_0=0$.

We consider $\mathbf{w}=\{w_0,w_1, \ldots\}$ as an adversarial disturbance with bounded energy. Let $\ell^k_2$ be the set of square-summable (bounded energy) signals of dimension $k$, i.e.,
$$
\ell^k_2 \!=\! \left\{\! \mathbf{w}\!=\!\{w_0,w_1, \ldots\} \mid \| \mathbf{w}\|_2 := \sqrt{\left(\textstyle \sum_{t=0}^\infty w_t^\tr w_t\right)} < \infty \right\}. 
$$
Robust $\mathcal{H}_\infty$ control aims to select the control sequence $\mathbf{u}=\{u_0, u_1, \ldots\}$ to minimize the quadratic performance $\sum_{t=0}^\infty (x_t^\tr Qx_t + u_t^\tr R u_t)$ against the worst disturbance of bounded energy $\|\mathbf{w}\|_2\leq 1$. Without loss of generality, we consider the energy bound on $\mathbf{w}$ to be 1, and we can formulate robust $\mathcal{H}_\infty$ control with any $\mathbf{w} \in \ell^{n_w}_2$. We thus consider the following min-max robust control problem
\begin{equation} \label{eq:robust-control-time-domain}   \min_{\mathbf{u}}\,\max_{\|\mathbf{w}\|_2\leq 1}\, \sum_{t=0}^\infty (x_t^\tr Qx_t + u_t^\tr R u_t).
\end{equation}
Throughout this paper, we make the following assumption 
\begin{assumption} \label{assumption:stablizability}
    The system matrices $(A, B)$ are stabilizable, and
    $B_w$ and $C$ are row full-rank.
    The weight matrices $Q, R$ are positive definite, and the controller has access only to the output sequence $\mathbf{y}:=\{y_0, y_1, \ldots\}$. 
\end{assumption}

Under \Cref{assumption:stablizability}, if full state measurements are available (i.e., $C=I_{n_x}$ in \cref{eq:dynamic}), it is well known that the optimal solution of \cref{eq:robust-control-time-domain} is a static state-feedback policy $u_t = Kx_t$, where $K\in\mathbb{R}^{n_u\times n_x}$ is a constant gain matrix \cite{zhou1996robust}. The general output-feedback $\mathcal{H}_\infty$ problem, however, is substantially more challenging \cite{glover2005state}. In this work, we restrict attention to \textit{linear static output-feedback policies} of the form 
\begin{equation} \label{eq:static-policies}
    u_t = K y_t, \qquad t= 0,1, 2, \ldots,
\end{equation}
and study the following policy optimization problem:
\begin{equation} \label{eq:robust-policy-optimization}
    \begin{aligned}
        \min_{K \in \mathbb{R}^{n_u \times n_y}} \; \max_{\|\mathbf{w}\|_2\leq 1}\quad &\sum_{t=0}^\infty (x_t^\tr Qx_t + u_t^\tr R u_t) \\
        \text{subject to} \quad & \text{\cref{eq:dynamic,eq:static-policies}}, \,x_0 = 0.
    \end{aligned}
\end{equation}

Given any feedback gain $K\in\mathbb{R}^{n_u \times n_y}$, the closed-loop dynamics from \cref{eq:dynamic,eq:static-policies} becomes 
 $   x_t \!=\! (A + BKC)x_t + B_ww_t,\, x_0 = 0.$
We denote the set of stabilizing gains as 
\begin{equation}
    \mathcal{K}:= \{K\in\mathbb{R}^{n_u\times n_y}\mid \rho(A+BKC) < 1 \}, 
\end{equation}
where $\rho(\cdot)$ denotes the spectral radius.

\begin{assumption}\label{assumption:K_nonempty}
    The set $\mathcal{K}$ is non-empty. 
\end{assumption}

If $K\in\mathcal{K}$, then the state trajectory $\mathbf{x}=\{x_0,x_1,\ldots\}$ belongs to $\ell_2^{n_x}$ for any bounded disturbance $\mathbf{w}\in\ell_2^{n_w}$. 
Let a performance signal be $z_t = (Q + C^\tr K^\tr R K C)^{1/2}x_t$ and $\mathbf{z} = \{z_0,z_1,z_2,\ldots\}$. Then, we have 
\begin{equation*}
    \|\mathbf{z}\|_2^2 = \sum_{t=0}^\infty (x_t^\tr Qx_t + u_t^\tr R u_t) < \infty, \quad \forall \mathbf{x}\in\ell^{n_x}_2, \mathbf{w}\in\ell^{n_w}_2. 
\end{equation*}

For any $K\in\mathcal{K}$, the closed-loop system defined by \cref{eq:dynamic,eq:static-policies} can be viewed as a bounded linear operator mapping disturbances $\mathbf{w}\in\ell_2^{n_w}$ to performance signals $\mathbf{z}\in\ell_2^{n_x}$. We denote this linear operator as $\mathbb{T}_{zw}(K)$, and define the $\ell_2 \to \ell_2$ induced norm as 
\begin{equation} \label{eq:linear-operator-time-domain}
\|\mathbb{T}_{zw}(K)\|:= \max_{\|\mathbf{w}\|_2\neq 0} 
\|\mathbf{z}\|_2/\|\mathbf{w}\|_2
\qquad \forall K \in \mathcal{K}.
\end{equation}
Since $\mathbb{T}_{zw}(K)$ is linear, it is not difficult to see that problem \cref{eq:robust-policy-optimization} is equivalent to $\min_{K \in \mathcal{K}} \, \|\mathbb{T}_{zw}(K)\|^2$. 
This linear operator \cref{eq:linear-operator-time-domain} is expressed in the time domain, and it also admits a frequency-domain representation via transfer functions. Denote transfer function matrix $\mathbf{T}_{zw}(K,\omega)$
from $\mathbf{w}$ to $\mathbf{z}$ as
\begin{equation*}
    \bT_{zw}(K,\omega)
    = 
    \begin{bmatrix}
Q^{1/2}\\
R^{1/2}KC\\
\end{bmatrix}
\left(
e^{j\omega}I- (A+BKC)
\right)^{-1}B_w,
\end{equation*}
where $\omega\in[0,2\pi]$ is a frequency variable. Its $\mathcal{H}_\infty$ norm is 
\begin{equation} \label{eq:Hinf-norm}
\| \bT_{zw}(K)\|_{\Hinf}
    = \max_{\omega\in [0,2\pi]}
    \sigma_{\max}(
    \bT_{zw}(K,\omega)
    ), \, \forall K \in \mathcal{K}
\end{equation}
where $\sigma_{\max}(\cdot)$ is the maximum singular value. The $\ell_2 \to \ell_2$ induced norm coincides with the $\mathcal{H}_\infty$ norm  \cite[Th. 4.4]{zhou1996robust}, i.e.,  
$
\|\mathbb{T}_{zw}(K)\| \!=\! \|\bT_{zw}(K)\|_{\Hinf}, \, \forall K \in \mathcal{K}.
$

Problem \cref{eq:robust-policy-optimization} is thus often called $\mathcal{H}_\infty$ optimization with static feedback policies, which can be equivalently~written~as
\begin{equation} \label{eq:policy-optimization-main}
    J^{\star}=\min_{K \in \mathcal{K}} \; J(K),
\end{equation}
where $J(K) :=  \|\bT_{zw}(K)\|_{\Hinf}$ is defined in \cref{eq:Hinf-norm}. We remark that problem \cref{eq:robust-policy-optimization} minimizes $J(K)^2$, whereas \cref{eq:policy-optimization-main} minimizes $J(K)$. Although the objective values differ by a square, the two problems clearly share the same set of optimal solutions.

\subsection{Problem statement}\label{subsection:problem_statement}

Robust $\mathcal{H}_\infty$ control is a cornerstone of control theory and has been studied extensively since the 1980s; see e.g., the classical textbook \cite{zhou1996robust}. In the state-feedback setting, standard approaches solve \cref{eq:policy-optimization-main} either through Riccati-based iterative methods \cite{doyle1988state}
or
via convex reformulations with Lyapunov variables \cite{gahinet1994linear}. The static output-feedback is substantially more challenging, and no method is known to compute its globally optimal solution. Both Lyapunov- and Riccati-based approaches rely on explicit knowledge of the system model. 

More recently, several works have explored direct optimization over the policy space $\mathcal{K}$ using local search methods \cite{apkarian2006nonsmooth,burke2006hifoo,saeki2006fixed,apkarian2009proximity,guo2022global,guo2023complexity,tang2023global}. These approaches are typically more scalable, and some of them are better-suited to model-free settings \cite{hu2023toward,talebi2024policy}. 
However, direct policy optimization for $\mathcal{H}_\infty$ control in \cref{eq:policy-optimization-main} is inherently a \textit{nonsmooth} and \textit{nonconvex} problem, even in the state-feedback setting. Because of these difficulties, the early work \cite{apkarian2006nonsmooth,burke2006hifoo,saeki2006fixed,apkarian2009proximity} provided no convergence-rate guarantees, while more recent studies \cite{guo2022global,guo2023complexity,tang2023global} established global optimality in the state-feedback case but relied on substantially more intricate algorithmic frameworks, often involving advanced randomized techniques~\cite{duchi2012randomized,zhang2020complexity,burke2005robust}.  

The main challenge is that the $\mathcal{H}_\infty$ cost function is both nonconvex and nonsmooth, with the nonsmoothness stemming from the maximization of the largest singular value in \cref{eq:Hinf-norm}. While it is known that this function is \textit{subdifferentially regular} \cite{apkarian2006nonsmooth,guo2022global,guo2023complexity,tang2023global}, such regularity is too broad to yield efficient algorithmic guarantees. At the same time, the function is more structured than a generic nonsmooth, nonconvex objective: it is naturally a \textit{lower-$C^2$} function \cite{rockafellar2009variational}, as can be seen directly from its definition \cref{eq:Hinf-norm}.  

This observation raises two central questions in this paper:  
\begin{enumerate}
    \item What additional structure (beyond subdifferential regularity) can be identified in the nonsmooth and nonconvex landscape of $\mathcal{H}_\infty$ policy optimization?  
    \item Given this structure, what are the simplest algorithms that can provably handle the nonconvex, nonsmooth nature of the $\mathcal{H}_\infty$ policy optimization problem?  
\end{enumerate}

To address the first question, we will show that the $\mathcal{H}_\infty$ cost function admits a form of \textit{weak convexity}, and in the state-feedback case, even satisfies a \textit{weak Polyak--\L{}ojasiewicz (PL)} condition. To address the second, we will demonstrate that classical first-order methods, such as the \textit{simple subgradient method}, can achieve deterministic non-asymptotic convergence rates that match the best-known~complexity for weakly convex functions \cite{davis2019stochastic,liang2023proximal,liao2025proximal,liao2024error}, while~avoiding randomized or overly complex procedures as in~\cite{guo2022global,guo2023complexity}.

\section{Weak Convexity in $\Hinf$ Optimization}\label{section:landscape}

We first summarize some basic landscape properties in \cref{eq:policy-optimization-main}, then establish the notion of weak~convexity, and finally discuss 
fundamental properties of
its stationary points.

\subsection{Basic landscape properties}

We here summarize the basic landscape properties~of~\cref{eq:policy-optimization-main}.

\begin{figure*}
\begin{subfigure}{.24\textwidth}
  \centering
\includegraphics[width=0.8\columnwidth]{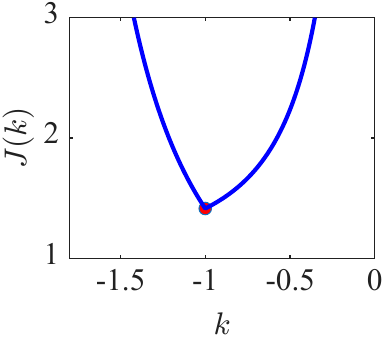}
\caption{$J(K)$ for system \Cref{eq:example-1}}
    \label{fig:1d}
\end{subfigure}%
\begin{subfigure}
{.24\textwidth}
    \centering
\includegraphics[width=0.95\linewidth]{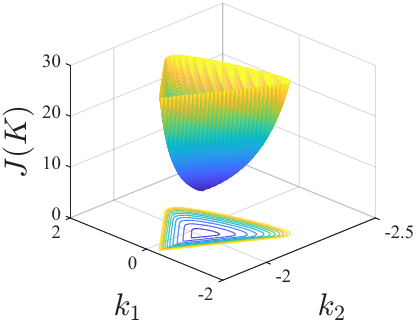}
    \caption{
$J(K)$ for system \Cref{eq:example-2}}
    \label{fig:desouza_xie}
\end{subfigure}
\begin{subfigure}{.24\textwidth}
  \centering
  \includegraphics[width=0.8\columnwidth]{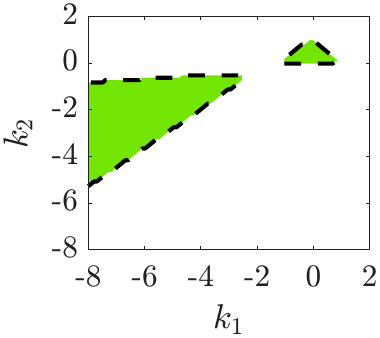}
\caption{
$\mathcal{K}$ in \cref{example:hinf_nonconvexity}}
    \label{fig:outputFB_K}
\end{subfigure}
\begin{subfigure}{.24\textwidth}
  \centering
\includegraphics[width=0.95\columnwidth]{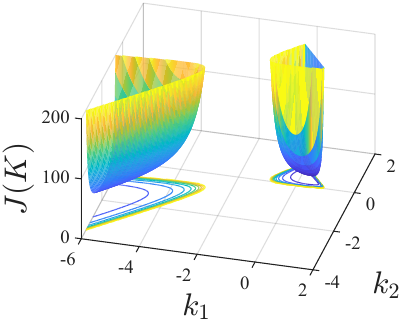}
\caption{
$J(K)$ in \cref{example:hinf_nonconvexity}}
    \label{fig:outputFB_J}
\end{subfigure}%
\caption{
Nonconvex and nonsmooth landscape in $\Hinf$ optimization.
(a)--(b) Nonsmoothness
of the cost functions $J$ in \cref{example:nonsmoothness};
(c)
Disconnectivity of 
$\mathcal{K}$
in \cref{example:hinf_nonconvexity}
with $\alpha=0.13$;
(d)
Nonconvexity of $J$
in \cref{example:hinf_nonconvexity}
with $\alpha=0.13$.
\vspace{-1mm}
}
\label{fig:hinf_plot}
\end{figure*}

\begin{lemma} \label{lemma:basic-properties}
    Suppose \cref{assumption:stablizability,assumption:K_nonempty} hold. 
    Consider the $\mathcal{H}_\infty$ optimization problem with static linear policies in \cref{eq:policy-optimization-main}. 
    Then the following statements hold:
    \begin{enumerate}
        \item The stabilizing set $\mathcal{K}$ is open and, in general, nonconvex. 
        If full state information is available (i.e., $C=I_{n_x}$), the set $\mathcal{K}$ is always path-connected, but in the general output-feedback case, it can be disconnected. 
        \item The $\mathcal{H}_\infty$ cost function $J:\mathcal{K}\to\mathbb{R}_+$ is coercive: 
        \[
        J(K)\to\infty \quad \text{as } K \to \partial\mathcal{K} \text{ or } \|K\|_F \to \infty.
        \]
        \item The function $J$ is generally nonsmooth and nonconvex. 
        \item On any sublevel set $\mathcal{K}_\nu := \{K \in \mathcal{K} \mid J(K) \leq \nu\}$, $J$ is $\beta_\nu$-Lipschitz continuous for some constant $\beta_\nu > 0$. 
    \end{enumerate}
\end{lemma}

These landscape properties are known in the literature; see, e.g., \cite{apkarian2006nonsmooth,guo2022global,guo2023complexity,tang2023global,zheng2023benign,feng2020connectivity}. We do not reproduce the proofs here, but instead highlight several insights and illustrative examples.  

First, it is straightforward to construct examples showing that $\mathcal{K}$ is nonconvex (as \cref{example:hinf_nonconvexity,example:saddle_local_min} below). When $C=I_{n_x}$, path-connectivity of $\mathcal{K}$ follows from an equivalent convex reformulation, since convex sets are always path-connected. 
Second, the coercivity of $J$ can be established by deriving a global quadratic lower bound; see \cite[Lem. 3.2]{guo2023complexity} for details. This property is crucial because it ensures that each sublevel set  
$
\mathcal{K}_\nu := \{K \in \mathcal{K} \mid J(K) \leq \nu\}
$ 
is compact.  Third, $J$ is naturally nonconvex since its domain $\mathcal{K}$ is already nonconvex. Its nonsmoothness arises from two sources in the definition \cref{eq:Hinf-norm}: (i) the maximum singular value of a complex matrix, and (ii) the maximization over the frequency variable.  Finally, while the local Lipschitz continuity of $J$ is well known, identifying a uniform Lipschitz constant on a sublevel set $\mathcal{K}_\nu$ is more subtle. Previous works such as \cite[Th.~5]{guo2022global} and \cite[Sec.~3]{guo2023complexity} assume this $\beta_\nu$-Lipschitz property directly. In fact, it can be rigorously established by explicitly characterizing the Clarke subdifferential of $J$.    

We next present two examples to illustrate \cref{lemma:basic-properties}.

\begin{example}[Nonsmoothness of function $J$]\label{example:nonsmoothness}
Consider 
\begin{equation} \label{eq:example-1}
 x_{k+1}  = x_k + u_k +w_k,\ y_k = x_k
 \end{equation}
with $u_k = kx_k$ and $Q=R=1$.
It is clear that $\mathcal{K} = \{k\mid |1+k|<1 \}
=(-2,0).$
It is easy to see 
\begin{align*}
    &{\bT_{zw}(k,\omega)^\her \bT_{zw}(k,\omega)}
    ={ \frac{{1+k^2}}{1+(1+k)^2-2(1+k)\cos\omega}},
\end{align*}
and thus the $\Hinf$ cost can be obtained explicitly as
\begin{equation*}
    J(k) = 
    \begin{cases}
    \frac{\sqrt{1+k^2}}{|k|},&k\in(0,-1]\\
    \frac{\sqrt{1+k^2}}{k+2},&k\in(-1,-2). 
    \end{cases}    
\end{equation*}
The function $J$ is clearly nonsmooth at $k=-1$, as illustrated in \cref{fig:1d}.  
Moreover, $J(k)\to\infty$ as $k \to 0$ or $k \to -2$, which demonstrates the coercivity property.

Consider another state-feedback system
\begin{equation} \label{eq:example-2}
    A=\begin{bmatrix}
        0 & 2\\
        4 & 0.2
    \end{bmatrix},\quad
    B = \begin{bmatrix}
        1\\
        0
    \end{bmatrix},
    \quad
    B_w=I_2
\end{equation}
with 
$Q = \diag(1,10^{-3})$ and
$R=1$.
From the classical $\mathcal{H}_\infty$ theory, we can compute the optimal value of $J(K)$ to be $J^*\approx8.327$.
We plot the function $J$ 
for $K=[k_1,k_2]$
in \cref{fig:desouza_xie}.
We can see nonsmooth points on the landscape.
\end{example}

\begin{example}[Nonconvexity of domain $\mathcal{K}$ and function $J$]\label{example:hinf_nonconvexity}
    Consider the following problem data:
\begin{align*}
    &A=
    \begin{bmatrix}
       1 -\alpha & -0.1 & -0.1 \\
0.1 & 1& 0 \\
0 & 0.1 & 1
    \end{bmatrix}
,\quad B=
\begin{bmatrix}
   0.1 \\
0 \\
0
\end{bmatrix},\\
&B_w = I_3,
\, C=
\begin{bmatrix}
    0 & 1 & 1 \\
1 & -1 & 1
\end{bmatrix},\,
Q = 0.01I_3,\, R = 0.01.
\end{align*}
For this system, the static output-feedback gain is of the form $K = [k_1,k_2]$.
When $\alpha  \approx 0.05\sim 0.13 $,
the feasible region $\mathcal{K}$ is disconnected,
which directly leads to the nonconvex landscape.
We plot set $\mathcal{K}$ and the function $J$ 
for $\alpha=0.13$
in \cref{fig:outputFB_K} and \cref{fig:outputFB_J}, respectively.
\end{example}

\begin{remark}[$\mathcal{H}_\infty$ optimization in continuous-time systems]\label{remark:continuous-time}
The continuous-time counterpart of \cref{eq:policy-optimization-main} is likewise nonsmooth and nonconvex. 
The set of stabilizing gains remains nonconvex, is path-connected when $C=I_{n_x}$, and can be path-disconnected in general. 
However, unlike the discrete-time case, the coercivity of $J$ does not necessarily hold in continuous time (see \cite[Example~2.3]{zheng2024benign}). 
As a result, its sublevel sets may fail to be compact, and it is unclear whether $J$ is $\beta_\nu$-Lipschitz continuous on such sets. \hfill $\square$
\end{remark}

\subsection{Lower-$C^2$ and weak convexity}

We now establish a key technical result: although the $\mathcal{H}_\infty$ cost $J$ is nonconvex, it is in fact \textit{weakly convex}. A precise statement will be given later in \Cref{theorem:local_wc}. 
Recall that a
function $f: \mathbb{R}^n \to \mathbb{R}$ is $\rho$-weakly convex if the function $x \mapsto f(x) + \frac{\rho}{2}\|x\|^2$ is convex \cite{davis2019stochastic,liao2024error}.~Weak convexity provides a broad framework that subsumes convex, smooth, and many composite nonconvex functions, and enables principled first-order algorithm design~\cite{davis2019stochastic,liao2024error,liao2025proximal,liang2023proximal}. 

Another related notion is that of \textit{lower-$C^2$} functions. Roughly speaking, a lower-$C^2$ function may be nonsmooth, but its nonsmoothness arises solely from taking a supremum of smooth $C^2$ components, making it significantly more structured than generic nonsmooth functions. 
Note that a function $f:\mathcal{D}\to\mathbb{R}$ is of class $C^2$ if it is twice differentiable on its open domain $\mathcal{D}$ and both $\nabla f$ and $\nabla^2 f$ are continuous.  The formal definition of lower-$C^2$ functions is as follows.

\begin{definition}[Lower-$C^2$ \cite{rockafellar2009variational}] \label{definition:lower-c2}
    We say a function $f:\mathcal{D}\to\RR$, defined on an open domain  $\mathcal{D} \subset \RR^n$, 
is lower-$C^2$, if 
for each $\bar{x}\in\mathcal{D}$, there exist a neighborhood $V$ of $\bar{x}$, a compact set $\mathcal{T}$, and a mapping 
$\varphi:\mathcal{T}\times V \to \mathbb{R}$ such that  
\begin{equation} \label{eq:lower-c2-definition}
    f (x) = \max_{t\in \mathcal{T}} \varphi(t,x),\quad
    \forall x\in V
\end{equation}
where, for each $t\in \mathcal{T}$, the function $\varphi(t,\cdot)$ is $C^2$ in $x$, 
and the mappings $\varphi(\cdot,\cdot)$, $\nabla_x \varphi(\cdot,\cdot)$, 
and $\nabla_x^2 \varphi(\cdot,\cdot)$ 
are jointly continuous on $(t,x)\in \mathcal{T}\times V$. 
\end{definition}

Note that the smooth functions $\varphi(t,\cdot)$ in \cref{eq:lower-c2-definition} may depend on the point $\bar{x}$. 
Thus, a lower-$C^2$ function can be \emph{locally} represented as the maximum of a family of smooth~$C^2$~functions.  
We will show that the $\mathcal{H}_\infty$ cost $J$ is lower-$C^2$. 

Before doing so, let us recall some simple examples. First, any quadratic function is lower-$C^2$, and so is the maximum of finitely many quadratics\footnote{
In this case, we regard $\mathcal{T}$
as a subset of a discrete topological space.}, e.g.,  $
f(x) = \max_{t\in \{1,\ldots,p\}} f_t(x)$. 
Second, the maximum eigenvalue function of symmetric matrices $\lambda_{\max}(\cdot)$ is also lower-$C^2$,~as
\begin{equation} \label{eq:example-eigenvalue}
\lambda_{\max}(X) = \max_{u^\tr u = 1} u^\tr X u,
\end{equation}
where each component $\varphi(u,X) = u^\top X u$ is linear in $X$ and hence $C^2$.  
Finally, the maximum singular value of a real matrix $Y \in \mathbb{R}^{p \times m}$ is lower-$C^2$, because
\begin{equation} \label{eq:example-signular-value}
\sigma_{\max}(Y) = \max_{\substack{u \in \mathbb{C}^{m},\, v \in \mathbb{C}^p \\ \|u\|=\|v\|=1}} \Re({v^\her}  Y u),
\end{equation}
where 
where $\Re(\cdot)$ represents the real part and,
each component $\varphi(u,v,Y) = \Re({u^\her} Y v)$ is linear (hence $C^2$) in $Y$.
These examples illustrate that many spectral functions naturally fall into the lower-$C^2$ class, which enables our analysis of the $\mathcal{H}_\infty$ cost below; recall its definition \cref{eq:Hinf-norm}.

Formally, we have the following result. 

\begin{lemma}\label{fact:lower_C2}
Suppose \cref{assumption:stablizability,assumption:K_nonempty} hold. 
Consider the $\mathcal{H}_\infty$ cost function $J:\mathcal{K}\to\mathbb{R}_+$ defined in \cref{eq:policy-optimization-main}. 
Then the following properties hold:
\begin{enumerate}
    \item $J$ is a lower-$C^2$ function. 
    \item $J$ is \emph{locally weakly convex}: for every $K_0 \in \mathcal{K}$, there exist constants $m>0$ and $r>0$ such that 
    \begin{equation} \label{eq:local-weak-convexity}
        K \;\mapsto\; J(K) + \tfrac{m}{2}\|K\|_F^2
    \end{equation}
    is convex over 
    \(\mathbb{B}_r(K_0) := \{K \mid \|K - K_0\|_F \leq r\}\).
    \item $J$ is \emph{subdifferentially regular}, i.e., its directional derivative coincides with its Clarke directional derivative at every point. 
\end{enumerate}
\end{lemma}

\begin{proof}
Similar to \cref{eq:example-signular-value}, we can rewrite the $\mathcal{H}_\infty$ cost as 
\begin{equation*}
J(K)=\max _{\omega \in [0,2\pi]} \max_{\substack{u \in \mathbb{C}^{n_w},\, v \in \mathbb{C}^{n_x+n_u} \\ \|u\|=\|v\|=1}} \Re\left[ v^\her \bT_{zw}(K,\omega)u\right].
\end{equation*} 
Define
$
\varphi(\omega,u,v,K) := \Re\!\left[ v^\her \bT_{zw}(K,\omega) u \right].
$ 
It is easy to verify that $\varphi$ is analytic in $K$ over the open set $\mathcal{K}$, and jointly continuous in $(\omega,u,v,K)$ on 
$[0,2\pi] \times \{u \in \mathbb{C}^{n_w} \mid\! \|u\|=1\} \times \{v \in \mathbb{C}^{n_x+n_u} \mid\! \|v\|=1\} \times \mathcal{K}$.  
Hence, we may represent~$J$~as
\begin{equation} \label{eq:hinf-smooth-componenet}
J(K) = \max_{(\omega,u,v) \in \Theta} \varphi(\omega,u,v,K),
\end{equation}
where
$
\Theta := [0,2\pi] \times \{u \in \mathbb{C}^{n_w} : \|u\|=1\} \times \{v \in \mathbb{C}^{n_x+n_u} : \|v\|=1\}
$
is compact.  
By \Cref{definition:lower-c2}, it follows that $J$ is lower-$C^2$.  
The local weak convexity 2) and subdifferential regularity 3) follow from 
\cite[Th.~10.31 and~Th. 10.33]{rockafellar2009variational}.
\end{proof}

The proof above is inspired by the continuous-time case \cite[Lem. 9]{noll2013bundle}. Note that the local weak convexity \cref{eq:local-weak-convexity} is a natural property of lower $C^2$ functions \cite[Th.~10.33]{rockafellar2009variational}. However, the weak convexity constant $m$ may depend on the point and its neighborhood, which may be hard to estimate. 

It is worth noting that the smooth components in \cref{eq:example-eigenvalue,eq:example-signular-value,eq:hinf-smooth-componenet} are \emph{global}, in the sense that they do not vary with the choice of reference point $K$. This structure is considerably stronger than general lower-$C^2$ functions in  \Cref{definition:lower-c2}. Consequently, we may anticipate better-behaved properties for the $\mathcal{H}_\infty$ cost than those established in \Cref{fact:lower_C2}. This is indeed the case: we can establish a uniform weak convexity constant over any convex subset in a sublevel set of $J$. 

\begin{theorem}[Weak convexity]\label{theorem:local_wc}
Suppose \cref{assumption:stablizability,assumption:K_nonempty} hold. Consider the $\mathcal{H}_\infty$ cost function $J:\mathcal{K}\to\mathbb{R}_+$ defined in \cref{eq:policy-optimization-main}. Let $\nu > 0$ and define a sublevel set $\mathcal{K}_\nu := \{K \in \mathcal{K} \mid J(K) \leq \nu\}\neq \emptyset$. Then, for any nonempty convex subset $V \subset \mathcal{K}_\nu$, there exists a constant $m_\nu>0$  such that the~function 
\[
K \;\mapsto\; J(K) + \tfrac{m_\nu}{2}\|K\|_F^2
\]
is convex over $V$. 
\end{theorem}

The proof is not difficult by recalling the fact that the composition of a convex function with an $L$-smooth function is weakly convex \cite{davis2019stochastic}. 
In our case, the maximum singular value function $\sigma_{\max}(\cdot)$ is convex on $\mathbb{R}^{ p\times n}$. We need to handle the complex transfer function and the maximization over the frequency variable $\omega$ in \cref{eq:Hinf-norm}. For this, we first use a well-known result from complex SDPs \cite[Exercise 4.42]{boyd2004convex} to equivalently eliminate complex numbers in $f(\cdot,\omega):=\sigma_\mathrm{max}(\bT_{zw}(\cdot,\omega))$ for each fixed $\omega \in [0,2\pi]$. Then, since $\mathcal{K}_\nu$ is compact thanks to the coercivity, over any convex subset $V \subset \mathcal{K}_\nu$, we can view $f(\cdot,\omega)$ as a composition of a convex function with a real-valued $L_\nu$-smooth function with some $L_\nu>0$, which is thus $m_\nu$-weakly convex with
some $m_\nu>0$ independent of $\omega$. Therefore, we have $J(K)+\frac{m_\nu}{2}\|K\|_F^2 = \max_{\omega\in[0,2\pi]}\{
f(K,\omega)+\frac{m_\nu}{2}\|K\|_F^2
\}$ is convex on any convex subset of $\mathcal{K}_\nu$. A detailed proof is given in \Cref{subsection:proof_wcvx}. 

\begin{example}[Weak convexity]\label{example:weak_convexity}
Consider 
$f:(-1,\infty)\to\mathbb{R}_+$:
\begin{equation*}
    f(x)
    = \sigma_\mathrm{max}\left(
    c(x)
    \right),\quad
    c(x)=\begin{bmatrix}
        1 & x/2\\
        (1+x)^{-1} & -2\cos
        x
    \end{bmatrix}.
\end{equation*}
This function is slightly different from the $\mathcal{H}_\infty$ cost \cref{eq:Hinf-norm}, as it has no frequency variable. However, it captures the essential composition of $\sigma_{\max}$ with a smooth function. Indeed, on any bounded interval $V\subset(-1,\infty)$, the mapping $c(\cdot)$ is $L_V$-smooth for some $L_V>0$, and thus $f(\cdot)$ is $m_V$-weakly convex on $V$ with some $m_V>0$.
To illustrate,
we plot $f(x)$ 
and $f(x)+\frac{5}{4}\|x-1\|^2$ 
on $V=[0,5]$ in \cref{fig:wcxv}.
Clearly, $f(x)$ is convexified by adding a quadratic perturbation  
$\frac{5}{4}\|x-1\|^2$.
\begin{figure}
    \centering    \includegraphics[width=0.8\linewidth]{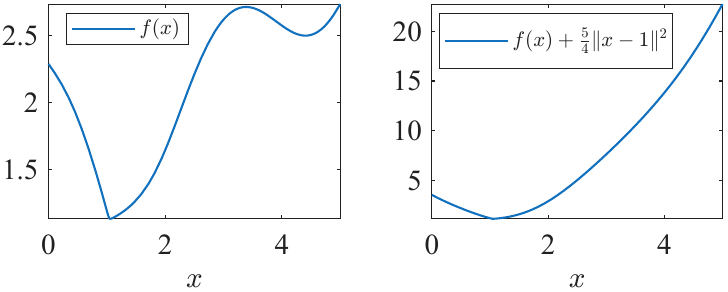}
\caption{
Plots of 
$f$ 
and $f+\frac{5}{4}\|\cdot-1\|^2$ 
on $V=[0,5]$
in
\cref{example:weak_convexity}.}
    \label{fig:wcxv}
\end{figure}
\end{example}

\subsection{Stationary points and weak PL property}

We next analyze the stationary points of $\Hinf$ cost $J$. 
As is typical in nonconvex optimization, spurious stationary points such as local minima and saddle points may arise in the general output-feedback $\mathcal{H}_\infty$ setting.  

\begin{lemma}
Suppose \cref{assumption:stablizability,assumption:K_nonempty} hold. 
Consider the $\mathcal{H}_\infty$ cost function $J:\mathcal{K}\to\mathbb{R}_+$ defined in \cref{eq:policy-optimization-main}.  When $C \neq I_{n_x}$, the function $J$ may admit spurious stationary points, including local minima and saddle points.   
\end{lemma}

We verify this result using \cref{example:saddle_local_min} below. This is
the same instance as \cref{example:hinf_nonconvexity} with a different choice of $\alpha$.
The example highlights that such spurious stationary points naturally emerge from the nonconvexity and potential disconnectedness of $\mathcal{K}$.

\begin{figure}
\begin{subfigure}{0.46\columnwidth}
  \centering
\includegraphics[width=0.85\columnwidth]{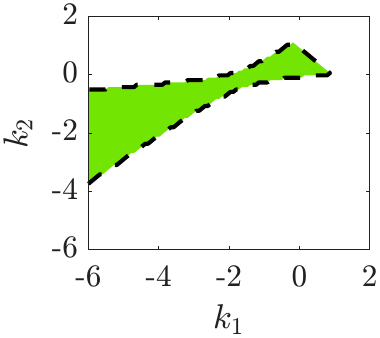}
\caption{$\mathcal{K}$
in \cref{example:saddle_local_min}
}
\label{fig:saddle_K}
\end{subfigure}%
\begin{subfigure}{0.52\columnwidth}
  \centering
\includegraphics[width=0.95\columnwidth]{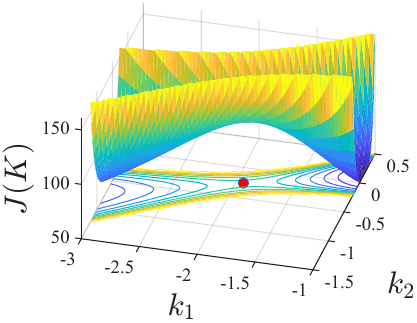}
\caption{
$J(K)$ in \cref{example:saddle_local_min}}
\label{fig:saddle_J_zoom}
\end{subfigure}
\caption{
Spurious stationary points of  $J$ 
with the set $\mathcal{K}$
in the case of static output-feedback in \cref{example:saddle_local_min}
(the same system as \cref{example:hinf_nonconvexity} with $\alpha=0.14$).
It can be observed that
$J(K)$ possesses not only a local minimum but also a saddle point, represented by the red dot.
}
\label{fig:saddle_point}
\end{figure}

\begin{example}[Saddle points]\label{example:saddle_local_min}
Consider the problem data 
in \cref{example:hinf_nonconvexity}.
For $\alpha=0.13$, the set $\mathcal{K}$ consists of two disconnected components (see \cref{fig:outputFB_K}). 
At $\alpha=0.14$, these components merge, and $\mathcal{K}$ becomes connected, as shown in \cref{fig:saddle_K}. 
In this case, we observe a saddle point in $J$ (see \cref{fig:saddle_J_zoom}), as well as at least one local minimum.  
\end{example}

In the full-state measurement case ($C=I_{n_x}$), the $\mathcal{H}_\infty$ optimization exhibits 
a much more tractable landscape, with a connected feasible region and convex-like shape as \cref{fig:1d,fig:desouza_xie}.
Indeed, we can establish a weak PL condition that ensures the global optimality of stationary~points. 

\begin{theorem}[Weak PL condition]\label{theorem:weak-PL}
Suppose \cref{assumption:stablizability,assumption:K_nonempty} hold.
If we have full state measurement, i.e., $C=I_{n_x}$, 
then, for $K\in\mathcal{K}$, there exists a constant
$\mu_K>0$
such that
\begin{equation}\label{eq:weakPL}
    \mu_K\left(J(K)-J^{\star}\right)
    \leq \mathrm{dist}(0,\partial J(K) ),
\end{equation}
where $\partial J$ denotes the Clarke subdifferential of $J$ and $J^\star = \min_{K \in \mathcal{K}} J(K)$. 
\end{theorem}

The proof is provided in \cref{subsubsection:proof_wcvx_wpl}.
We also provide an explicit form of the constant $\mu_K$ in \cref{eq:weakPL} in the proof. 
This remarkable property has a tight connection to the well-known LMI reformulation via the bounded real lemma \cite{boyd1994linear},
which will be utilized in the proof. Our proof is motivated by the recent ECL framework \cite{zheng2023benign} and the result in \cite{umenberger2022globally}.
The global optimality of stationary points follows as a corollary.
\begin{corollary}[Global optimality of stationary points]\label{corollary:stationary_global_optimality}
Under the same assumptions as \cref{theorem:weak-PL}, we have
$0\in \partial J(K)
\Leftrightarrow
    K\in {\arg\min}_{K\in\mathcal{K}}J(K).$
\end{corollary}

\begin{proof}
If $0 \in \partial J(K)$, then we have $\mathrm{dist}(0,\partial J(K) ) = 0$, which implies that $J(K) = J^\star$ by the weak PL \cref{eq:weakPL}. On the other hand, if $K\in {\arg\min}_{K\in\mathcal{K}}J(K)$, we must have $K$ to be stationary, and thus $0 \in \partial J(K)$ \cite[Th. 10.1]{rockafellar2009variational}.
\end{proof}

\subsection{Proof
of \cref{theorem:local_wc}}\label{subsection:proof_wcvx}

To establish \cref{theorem:local_wc}, we first rewrite the $\Hinf$ cost as the composition of  $\sigma_\mathrm{max}(\cdot)$ and
a real-valued mapping using \cref{lemma:complex-SDP} below,
which follows from a famous fact
\cite[Exercise 4.42]{boyd2004convex} in complex SDPs. 
We thereby identify the weakly convex constant $m_\nu>0$ via the compactness of $\mathcal{K}_\nu$
from \cref{lemma:basic-properties}.
The proof of \cref{lemma:complex-SDP} is presented in \cref{appendix:proof_complex_sdp} due to the page limit.
\begin{lemma}\label{lemma:complex-SDP}
Let $X$ be a $q\times r$ hermitian matrix. 
We have
\begin{equation*}
\sigma^2I \succeq 
X^\her X \quad\Longleftrightarrow\quad
\sigma I\succeq 
\Psi(X)
\end{equation*}
where $\Psi(\cdot)$ is a real-valued linear mapping
defined by
$$\Psi (X) = 
\small
\begin{bmatrix}
        \begin{bmatrix}
            0 & \FR(X)^\tr \\
            \FR(X)&0
        \end{bmatrix} & \begin{bmatrix}
            0 & \FI(X)^\tr \\
            -\FI(X)&0
        \end{bmatrix}\\
        \begin{bmatrix}
            0 & -\FI(X)^\tr  \\
            \FI(X)&0
        \end{bmatrix}
        &\begin{bmatrix}
            0 & \FR(X)^\tr \\
            \FR(X)&0
        \end{bmatrix} 
    \end{bmatrix},$$
where $\FR(X)$ and $\FJ(X)$ represent the real and imaginary parts, respectively.
\end{lemma}

Recall the definition of the function $J(K)=\|\bT_{zw}(K)\|_{\Hinf}$ in \cref{eq:Hinf-norm}.
We first show the weak convexity of
$$f(\cdot,\omega) := \sigma_{\max}(\bT_{zw}(\cdot,\omega))$$ for any fixed $\omega\in[0,2\pi]$.
Then, the weak convexity of $J = \max_{\omega \in [0,2\pi]} f(\cdot,\omega) $ follows from the standard result
for point-wise maximum of convex functions \cite[Sec. 3.2.3]{boyd2004convex}.

Fix an $\omega\in [0,2\pi]$.
By applying \cref{lemma:complex-SDP},
we have that
\begin{equation*}
f(K,\omega) = \sigma_{\max}(\bT_{zw}(K,\omega))
    =\lambda_{\max}
    \left(
    (\Psi\circ \bT_{zw})(K,\omega)
    \right).
\end{equation*}
It is clear that both $\FR(\bT_{zw}(K,\omega))$
and $\FI(\bT_{zw}(K,\omega))$
are real analytic on $\mathcal{K}\times[0,2\pi]$.
Thus, the mapping
$\Psi\circ \bT_{zw}$ is also real analytical over $\mathcal{K} \times [0,2\pi]$.
Then, the Hessian of $\Psi\circ \bT_{zw}$ exists and is continuous 
w.r.t. $(K,\omega)$
over $\mathcal{K}\times [0,2\pi]$.

By the compactness of $\mathcal{K}_\nu$,
we have the $L_\nu$-smoothness
of $f(\cdot,\omega)$ on $\mathcal{K}_\nu$ 
with the following finite positive constant:
\begin{align*}
L_{\nu}
= 
4\eta^2
\max_{(i,j)}
\max_{
\substack{
K\in\mathcal{K}_\nu,\\
\omega\in[0,2\pi]
}}
\left\|\nabla^2_{\mathrm{vec}(K)} [(\Psi\circ \bT_{zw})]_{ij}(K,\omega)\right\|_F,
\end{align*}
where $\eta = n_x+n_u+n_w$.
By recalling that $\lambda_{\max}(\cdot)$ is convex and $1$-Lipschitz, $f(K,\omega)$ is $m_\nu$-weakly convex over any convex subset of $K_\nu$ with
$m_\nu=L_\nu$; see \cite[Sec. 2.1]{davis2019stochastic}.

Finally, it is straightforward to see that the function 
\begin{equation*}
    J(K)+ \frac{m_\nu}{2}\|K\|_F^2
    =
    \max_{\omega\in[0,2\pi]}
    \left\{
    f(K,\omega) + \frac{m_\nu}{2}\|K\|_F^2
    \right\}
\end{equation*}
is convex on any convex subset of $\mathcal{K}_\nu$ from the $m_\nu$-weak convexity of each $f(K,\omega)$
\cite[Sec. 3.2.3]{boyd2004convex}. This leads to the desired $m_\nu$-weak convexity of $J(K)$. 
\hfill$\square$

\section{The Subgradient Method}\label{section:subgradient_method}

In this section, we discuss the subgradient method to solve the $\mathcal{H}_\infty$ optimization \cref{eq:policy-optimization-main}.
Thanks to the weak convexity in \Cref{theorem:local_wc}, we provide a non-asymptotic convergence rate under mild assumptions, by adapting the recent result~in~\cite{davis2019stochastic}. 

\subsection{The subgradient method}

Given an initial point $K_0 \in \mathcal{K}$, the subgradient method simply follows the update below 
\begin{equation}\label{eq:subGM}
    K_{t+1} = K_t - \alpha_{t} G_t,\quad
    G_t\in \partial J(K_t),
\end{equation}
where $\partial J(K)$ is the Clarke subdifferential of $J(\cdot)$ at $K\in \mathcal{K}$ \cite{rockafellar2009variational}.
This is arguably the simplest algorithm for nonsmooth nonconvex optimization.
Thanks to its simplicity, we may implement \eqref{eq:subGM} in a model-free manner relatively easily, only with a subgradient oracle.

\subsection{Non-asymptotic convergence guarantees}

Despite its simplicity, the subgradient method is classically known to achieve sublinear convergence only for convex and/or smooth problems \cite{nesterov2018lectures,davis2019stochastic,zhang2020complexity}. Its behavior in the nonsmooth nonconvex setting was far less understood
because (i) stationary concepts in convex and/or smooth cases can no longer be used
\cite{davis2019stochastic,zhang2020complexity},
and~(ii) the subgradient update \cref{eq:subGM} does not guarantee monotonic decrease of the cost; it is indeed possible even in the nonsmooth convex case that $J(K_{t+1}) > J(K_t)$ for~any $\alpha_t > 0$.
Recently, \cite{davis2019stochastic} established the first non-asymptotic rate for minimizing 
\textit{unconstrained} weakly convex functions. In this work, we extend their analysis to the $\mathcal{H}_\infty$ optimization \cref{eq:policy-optimization-main}
involving \emph{a nonconvex
stabilizing constraint} $\mathcal{K}$.

One key concept is the \textit{Moreau envelope}, which allows for quantifying the performance of \cref{eq:subGM} even in the absence of a function value decrease.
\begin{definition}[Moreau envelope]
For $K\in\mathcal{K}$ and
$\rho>0$,
the Moreau envelope $J_\rho(\cdot)$ is defined as
    $J_\rho(K) := \min_{M\in\mathcal{K}} J(M) + 
    \frac{\rho}{2} \|M-K\|_F^2$.
\end{definition}

By the coercivity in \cref{lemma:basic-properties},
there always exists a minimizer
$\hat{K}\in\arg\min_{M\in\mathcal{K}} J(M) + 
    \frac{\rho}{2} \|M-K\|_F^2$, but the minimizer may not be unique. By the definition, we have 
$$J_\rho(K)\leq J(K), \;  \text{and} \; J_\rho(\hat K)=J(\hat K)\leq J(K), \quad \forall K \in \mathcal{K}.$$
While the Moreau envelope is convenient in theory, 
computing the value of $J_\rho(\cdot)$ and $\hat{K}$ might be untractable in practice.
If $\rho>0$ is sufficiently large, the minimizer $\hat{K}$ will be unique. We have the following result by adapting  \cite[Sec. 2.2]{davis2019stochastic}.
\begin{lemma}\label{lemma:Moreau}
Consider the $\mathcal{H}_\infty$ cost function $J:\mathcal{K}\to\mathbb{R}_+$ defined in \cref{eq:policy-optimization-main}.  Fix any $K\in\mathcal{K}$. Let $\nu>J(K)$, and define $L>0$ and $m>0$ such that
$J(\cdot)$ is
$L$-Lipschitz over $\mathcal{K}_\nu$ and
$m$-weakly convex over any convex subset of $\mathcal{K}_\nu$.
Then, for
$\rho>
\max\{m,L/
\mathrm{dist}(K,\partial \mathcal{K}_\nu)
\}$,
we have
\begin{align*}
\|K - \hat{K}\|_F &\leq \mathrm{dist}(K,\partial \mathcal{K}_\nu) \\
   \nabla J_\rho( K)
   &
   = \rho(K-\hat K),\; 
   \\
\mathrm{dist}(0,\partial J(\hat K))
&
\leq \| \nabla J_\rho(K)\|_F,
\end{align*}
where $\{\hat{K}\}=\arg\min_{M\in\mathcal{K}} J(M)+\frac{\rho}{2}\|M-K\|_F^2$.
\end{lemma}
\begin{proof}
 From the optimality condition of $\min_{M} J(M) + \frac{\rho}{2}\|M - K\|_F^2$, we have 
 \begin{align}
    \label{eq:opt-cnd-env}
     \rho ( K - \hat{K}) \in \partial J(\hat{K}),
 \end{align}
 where $\hat{K}$ is any minimizer. The $L$-Lipschitzness of $J$ implies
     $\|K - \hat{K}\|_F \leq L /\rho \leq d:= \mathrm{dist}(K,\partial \mathcal{K}_\nu)$,
 where the last inequality is by the choice of $\rho >\max\{m,L/
\mathrm{dist}(K,\partial \mathcal{K}_\nu)\}.$ Thus, minimizing $J$ over $\mathcal{K}$ is the same as minimizing $J$ over the ball $\mathbb{B}_{d}(K)$, i.e.,
$ \arg\min_{M \in \mathcal{K}} \{ J(M) + \frac{\rho}{2}\|M - K\|_F^2\}
    = \arg\min_{M \in \mathbb{B}_{d}(K) }\{ J(M) + \frac{\rho}{2}\|M - K\|_F^2\}.$

Since $\mathbb{B}_{d}(K) \subseteq \mathcal{K}_{\nu}$ and $\rho >m$, we have $J(\cdot) + \frac{\rho}{2}\|\cdot - K\|_F^2$ being $(\rho - m)$-strongly convex over $\mathbb{B}_{d}(K)$. The minimizer is then unique. 
Hence, from \cite[Lem. 2.2]{davis2019stochastic}, the gradient of $J_{\rho}$ can be computed as 
   $\nabla J_{\rho}(K) =  \rho (K-
   \hat{K} )$,
and \cref{eq:opt-cnd-env} shows $\mathrm{dist}(0,\partial J (\hat{K}) ) \leq \| \nabla J_{\rho}(K)\|_F$.
\end{proof}

We are now ready to state the convergence of \cref{eq:subGM} in terms of the Moreau envelope
under a technical condition that the costs are bounded along the iterations. 

\begin{theorem}\label{theorem:subgradient_descent}
Consider the $\mathcal{H}_\infty$ cost function $J:\mathcal{K}\to\mathbb{R}_+$ defined in \cref{eq:policy-optimization-main}. Let 
 $\{K_0,\ldots,K_{T+1}\}$ be a sequence generated by
the subgradient method \cref{eq:subGM} with $K_0\in \mathcal{K}$
and step sizes
$\{\alpha_0,\ldots,\alpha_T\}$.
Suppose
$\max_{t\in\{0,\ldots,T+1\} }J(K_t) <\tilde{J}<\infty$
and $\max_{t\in{0,\ldots,T}}\|G_t\|_F \leq H<\infty.$
Then, we have
the feasibility of the iterates
$\{K_t\}\subset\mathcal{
K
}$, and
\begin{align*}
    \min_{t=0,\ldots,T}
    \| \nabla J_{\rho}(K_t)\|_F^2    
    \leq
   \kappa
    \frac{
    J_{\rho}(K_0)-J^\star
    + \frac{
    {\rho}_{\min}
    H^2}{2} \sum_{t=0}^T \alpha_k^2
    }{
    \sum_{t=0}^T \alpha_t},
\end{align*}
where
$\kappa = \rho/(\rho-m)$
with some constant 
$m>0$, and 
$\rho > {\rho}_{\min}\in(m,\infty)$.
Furthermore, if $\alpha_0=\cdots =\alpha_{T}= \alpha = \beta /\sqrt{T+1}$ with $\beta>0$,
we have
\begin{equation*}
   \min_{t=0,\ldots,T}
    \| \nabla J_{\rho}(K_t)\|_F^2  
    \leq 
    \kappa
    \frac{
    J_{\rho}(K_0)-J^\star
    +{\rho}_{\min} H^2\beta^2/{2}} 
    {
    \beta \sqrt{T+1}
    }.
\end{equation*}
\end{theorem}
\begin{proof}
Consider the sublevel set $\mathcal{K}_{\tilde{J}}
=
\{
K\in 
\mathcal{K} \mid J(K) \leq \tilde{J} \}$, which is compact due to the coercivity. By \Cref{theorem:local_wc,lemma:basic-properties}, 
$J(\cdot)$ is 
$L$-Lipschitz
over $\mathcal{K}_{\tilde{J}}$ 
and
$m$-weakly convex 
over any convex subset of $\mathcal{K}_{\tilde{J}}$, with some 
$L,\,m>0$.

Let ${\rho}_{\min}>
\max\{m,L/\tilde d\},
$
where  
$$\tilde d:=
\min_{t\in\{0,\ldots ,T+1\}} \mathrm{dist}(K_t,\partial \mathcal{K}_{\tilde J}) > 0.
$$
Then, 
the same argument as \cref{lemma:Moreau}
implies that
$\mathbb{B}_{\tilde{d}}(K_t)\subset \mathcal{K}_{\tilde{J}}$ for all $t$, and
 $\hat{K}_t\in \arg\min_{M}J(M)+\frac{\rho}{2}\|M-K_t\|_F^2$
is contained by the ball $\mathbb{B}_{\tilde{d}}(K_t)$.
Hence,
we can regard $J(K) + \frac{\rho}{2}\|K-K_t\|_F^2$ as 
a $(\rho-m)$-strongly convex function over
$\mathbb{B}_{\tilde d}({K}_t)$, and
this theorem follows from 
performing the same argument as \cite[Th. 3.1]{davis2019stochastic}.
The coercivity  in \cref{lemma:basic-properties} yields
the feasibility of $\{K_t\}$.
\end{proof}

This theorem indicates that
we can find an iterate $K_{t^*}$ satisfying
    $\|\nabla J_{\rho}(
    K_{t^*})\| \leq \epsilon$
after
$$
T= \!
\left\lceil 
\kappa^2\left(J_{\rho}(K_0)\!-\!J^\star
    \!+\!
    \textstyle\frac{m H^2\beta^2}{2} \right)^2
    \beta^{-2}\epsilon^{-4}
\!-1\right\rceil
\!=\mathcal{O}(\epsilon^{-4})
$$ iterations,
as long as the costs remain finite.
For $K_{t^*}$,
\begin{equation*}
    \mathrm{dist}(0,\partial J(\hat K_{t^*})) \leq \epsilon,\quad
    \|\hat K_{t^*}-K_{t^*}\|_F\leq \epsilon/\rho
\end{equation*}
hold,
where $\hat{K}_{t^*} = 
 {\arg\min}_{M\in \mathcal{K}} J(M)+\frac{\rho}{2}\|M-K_{t^*}\|_F^2
$.
Thus, 
the subgradient method \cref{eq:subGM} allows us to
arrive at an approximation of a stationary point 
in the above sense.

\begin{remark}
\cref{theorem:subgradient_descent} presents
a non-asymptotic convergence rate 
for the subgradient method \cref{eq:subGM}.
Notice that the boundedness of the cost is directly assumed, as it is often
empirically true 
for sufficiently small $\{\alpha_t\}$
but 
not formally guaranteed in \cref{eq:subGM}.
Under this assumption,
we can ensure the feasibility $\{K_t\}\subset\mathcal{K}$ by the coercivity in \cref{lemma:basic-properties},
bypassing the nonconvex constraint set $\mathcal{K}$
through a careful choice of the parameter $\rho$.
To alleviate the assumption on the cost, one may use other algorithms that can generate a descent direction (e.g., the proximal bundle method \cite{liao2025proximal}). We leave it as our future work. \hfill $\square$
\end{remark}

\section{Numerical results}\label{section:simulations}

Here, we run the subgradient method \cref{eq:subGM}
with fixed step sizes $\alpha_t = \alpha^1$ or $\alpha^2$,
where $(\alpha^1,\alpha^2)=(10^{-3},10^{-4})$,
for both full state-feedback and static output-feedback cases.
The following results demonstrate the effectiveness of the subgradient method \cref{eq:subGM} and validate our theoretical results.

\begin{figure}
\begin{subfigure}{0.48\columnwidth}
  \centering
\includegraphics[width=\columnwidth]{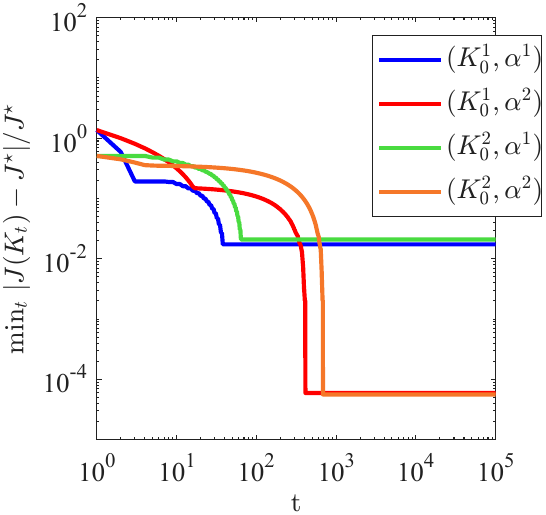}
\caption{
System \cref{eq:example-2}
($C= I_{n_x}$)
}
\label{fig:sm_smfb}
\end{subfigure}
\begin{subfigure}{0.48\columnwidth}
  \centering
\includegraphics[width=\columnwidth]{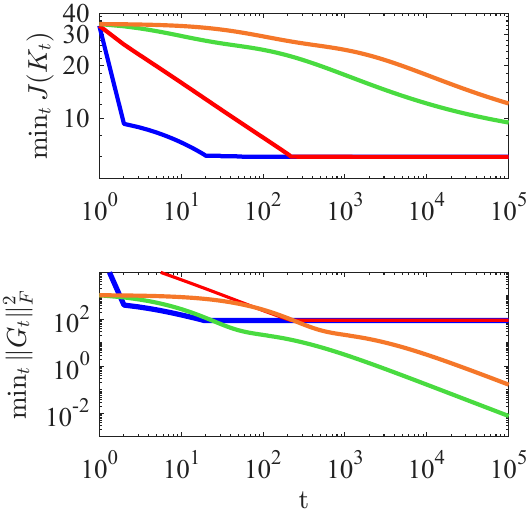}
\caption{
\cref{example:saddle_local_min} ($C\neq I_{n_x}$)
}
\label{fig:sm_ofb}
\end{subfigure}%
\caption{
The simulation results of the subgradient method \cref{eq:subGM}
with
$\alpha_t=\alpha^1=10^{-3}$ or $\alpha^2=10^{-4}$:
(a) The plot of 
$\min_t |J(K_t)-J^\star|/J^\star$
for the system \cref{eq:example-2} 
with a full state-feedback
in \cref{example:nonsmoothness}
 starting from $K_0=K_0^1$ or $K_0^2$,
 where $K_0^1=[0, -1.9]$ and $K_0^2=[0.2,-2]$;
(b) The 
plots of 
$\min_t J(K_t)$ (upper)
and
$\min_t \|G_t\|_F^2$ 
with $G_t\in\partial J(K_t)$ (lower)
for \cref{example:saddle_local_min} with $K_0 = K_0^1=[0, 0]$ or $K_0^2=[-5,-2]$.
}
\label{fig:sm_simulations}
\end{figure}

\textit{Full state-feedback.}
For the system \cref{eq:example-2} in \cref{example:nonsmoothness},
we test \cref{eq:subGM} with different
parameter choices.
Specifically, we use
$K_0 = K_0^1=[0, -1.9]$ or $K_0^2=[0.2,-2]$.
We plot the simulation result in \cref{fig:sm_smfb}, where the vertical axis shows the best relative objective residual $\min_t |J(K_t)-J^\star|/J^\star$ at $t$.
In all the cases, the subgradient method generates feasible iterates, which implies that the assumptions in \cref{theorem:subgradient_descent} are satisfied.
In addition,
as expected, we obtained a globally optimal solution 
for both initial conditions with favorable accuracy,
especially when $\alpha_t=\alpha^2$.

\textit{Static output-feedback.}
For  \cref{example:saddle_local_min},
we similarly run \cref{eq:subGM} with four 
different choices of the initial condition and step size.
In particular, we use
$K_0 = K_0^1$ or $K_0^2$
with $K_0^1=[0,0],\,K_0^2=[-5,-2]$.
The simulation result is presented in \cref{fig:sm_ofb}.
It can be seen that
the feasibility is satisfied in all the cases.
Regarding the performance,
the cost value converged to 
a smaller value
in the case of $K_0=K_0^1$.
In addition, when $K_0=K_0^2$,
we observe that the subgradient
also converged nearly to a stationary point. 
When $K_0=K_0^1$ in \cref{fig:sm_ofb},
the magnitude of the subgradient
remained larger despite the smaller cost value.
This is due to a zigzagging behavior around a stationary point and
may happen when it is nonsmooth.
For example, $f(x)=|x|$ takes the minimum at $x=0$, but the origin is nonsmooth and $\partial f(0)=[-1,1]\neq \{0\}$,
which may lead to such behavior.

\section{Conclusion}\label{section:conclusion}

This study addressed the $\Hinf$ optimization problem for discrete-time linear systems with static output-feedback.
We first proved the weak convexity
and further established 
the weak PL condition in the full state-feedback case.
Then, we also discussed the subgradient method from a weak convexity perspective, presenting the first non-asymptotic convergence rate. The numerical simulations validated our theoretical analysis and demonstrated the effectiveness of the subgradient method.
Our future work includes alleviating the assumption of the cost value boundedness in the convergence rate result, for example, through the proximal bundle method.

\bibliographystyle{IEEEtran}
\bibliography{ref.bib}

\newpage
\appendix

\subsection{Proof of \cref{lemma:complex-SDP}}\label{appendix:proof_complex_sdp}

Consider $X\neq 0$ without loss of generality.
By the Schur complement,
we have
\begin{equation}\label{eq:proof_complexSDP1}
    \sigma I \succeq \frac{1}{\sigma}X^\her X\;
\Longleftrightarrow\;
\sigma I -
\begin{bmatrix}
    0 & X^\her \\
    X & 0
\end{bmatrix}
\succeq 0.
\end{equation}
Now, 
from \cite[Exercise 4.42]{boyd2004convex}
for a hermitian matrix $Z$,
we know that $\lambda I \succeq  Z$ is equivalent to
$\lambda I\succeq  
   \begin{bmatrix}
       \mathfrak{R}(Z)  & -\mathfrak{I}(Z)\\
       \mathfrak{I}(Z)&
       \mathfrak{R}(Z)
   \end{bmatrix}$.
With this, letting $    Y =
    \begin{bmatrix}
    0 & X^\her \\
    X & 0
    \end{bmatrix}$ in \cref{eq:proof_complexSDP1}, 
we obtain
\begin{equation*}
\sigma I-Y
\succeq 0\quad\Longleftrightarrow
\quad
\sigma I\succeq 
\begin{bmatrix}
    \mathfrak{R}(Y)  & -\mathfrak{I}(Y) \\ 
    \mathfrak{I}(Y) & \mathfrak{R}(Y)  
\end{bmatrix}
= \Psi(X). 
\end{equation*}
This completes the proof. 
\hfill$\square$

\subsection{Proof of \cref{theorem:weak-PL}}\label{subsubsection:proof_wcvx_wpl}

We derive the weak PL condition \cref{eq:weakPL} by exploiting 
a favorable \textit{partial minimization} structure of the cost $J$.
First,
by using the nonstrict bounded real lemma \cite[Sec.  2.3.7]{boyd1994linear}, 
it is known that
the problem \cref{eq:policy-optimization-main} with $C=I_{n_x}$ admits the following reformulation with the additional Lyapunov variable $P$ and a scalar $\gamma$:
\begin{equation}\label{eq:Hinf_opt_lif}
    \min_{\gamma,K,X}\quad \gamma \quad
    \text{subject to}\quad
    (K,\gamma,P)\in \mathcal{L}_\mathrm{lft},
\end{equation}
where we define the set $\mathcal{L}_\mathrm{lft}$ as 
\begin{equation*}
    \mathcal{L}_\mathrm{lft}=
    \left\{
    (K,\gamma,P)
    \mid
    P\succ 0,\,
    \Lambda(K,\gamma,P)\preceq 0
    \right\}
\end{equation*}
with
\begin{align*}
\Lambda(K,\gamma,P)
\!=\!\begin{bmatrix}
    A_K^\tr PA_K-P & \!\!\! A_K^\tr P \\
PA_K &  \!\!\!P
\end{bmatrix}
\!+\!
\begin{bmatrix}
    Q+K^\tr R K & \!\!\!\!\!\! 0 \\
0 &  \!\!\!\!\!\! -\gamma^2 I
\end{bmatrix}
\end{align*}
where $A_K=A+BK$.
Then, we can represent $J(K)$ as the partial minimization of the following lifted function:
\begin{equation*}
    J_\mathrm{lft}(K,\gamma,P)
    =  
    \gamma + \delta_{\mathcal{L}_\mathrm{lft}}
    (K,\gamma,P)
\end{equation*}
with the indicator function $\delta_{\mathcal{L}_\mathrm{lft}}
    (K,\gamma,P)$ for $\mathcal{L}_\mathrm{lft}$.
Namely,
\begin{equation}\label{eq:Jk_parmin}
    J(K) = \inf_{\gamma
    ,P}
    J_\mathrm{lft}(K,\gamma,P).
\end{equation}

We first show the weak PL inequality for $J_\mathrm{lft}(K,\gamma,P)$,
by which we derive \cref{eq:weakPL}.
We here use the exact convex reformulation (LMI) \cite{boyd1994linear} for \cref{eq:Hinf_opt_lif}
with
the change of variables
$X=\left(P/\gamma\right)^{-1}$
and 
$Y=KX$.
To handle this, we define the following $C^\infty$ mapping $\Pi:\mathcal{L}_\mathrm{lft}\to \mathcal{F}_\mathrm{cvx}$:
\begin{equation*}
    \Pi(K,\gamma,P) = \left(\gamma,
    K(P/\gamma)^{-1},(P/\gamma)^{-1}
    \right),
\end{equation*}
with
the following convex set
\begin{equation*}
    \mathcal{F}_\mathrm{cvx}=
    \{
    (\gamma,Y,X)\mid
    X\succ 0,\,
    \mathrm{LMI}(\gamma,Y,X)\preceq 0
    \},
\end{equation*}
where $\mathrm{LMI}(\gamma,Y,X)$ is the following linear map:
\begin{align*}
    &\mathrm{LMI}(\gamma,Y,X)=\\
    &
    \begin{bmatrix}
-X & 0 & X &  \!\!\! (A X+B Y)^\tr & Y^\tr \\
0 & -\gamma I & 0 & \!\!\! I & 0 \\
X & 0 & -\gamma Q^{-1} & \!\!\! 0 & 0 \\
(A X+B Y) & I & 0 & \!\!\! -X & 0 \\
Y & 0 & 0 & \!\!\! 0 & -\gamma R^{-1}
    \end{bmatrix}.
\end{align*}
This LMI is standard and  follows from
the inequalities $\Lambda(K,\gamma,P)\preceq 0,\,P\succ 0$ and
Schur complement,
which yields
$$\Pi(\mathcal{L}_\mathrm{lft})
    = \mathcal{F}_\mathrm{cvx}.$$
Notice that
this mapping $\Pi(\cdot,\cdot,\cdot)$ is a diffeomorphism.
Then, by defining the function
\begin{equation*}
    J_\mathrm{cvx}(\gamma,Y,X)
    = 
    \gamma
    + \delta_{\mathcal{F}_\mathrm{cvx}}(\gamma,Y,X),
\end{equation*}
which is convex from the convexity of $\mathcal{F}_\mathrm{cvx}$,
we can further rewrite 
\begin{equation*}
    J_\mathrm{lft}(K,\gamma,P) = 
    (J_\mathrm{cvx}\circ \Pi)
    (K,\gamma,P).
\end{equation*}
Note also that
$J_\mathrm{cvx}(\gamma,Y,X) = 
    (J_\mathrm{lft}\circ \Pi^{-1})
    (\gamma,Y,X).$

Now we are ready to establish the weak PL condition.
By the convexity of $J_\mathrm{cvx}$,
\begin{align}
\label{eq:convexity_J_cvx}
\begin{aligned}
    &
    J_\mathrm{lft}(K,\gamma,P)-J^\star=
    (J_\mathrm{cvx}\circ \Pi)
    (K,\gamma,P)- J^\star
    \\
    \leq& 
    \langle W,
    \Pi(K,\gamma,P)  
    -\Pi(K^\star,J^\star,P^\star)
    \rangle\\
\leq&
\|W\|_F
    \times \|\Pi(K,\gamma,P)    -\Pi(K^\star,J^\star,P^\star)\|_F
    \end{aligned}
\end{align}
where 
$W\in \partial J_\mathrm{cvx}
    \left(
    \Pi(K,\gamma,P)
    \right)$, and
$(K^\star,J^\star,P^\star)$ represents an optimal triplet.
The chain rule 
\cite[Th. 10.6]{rockafellar2009variational}
for $J_\mathrm{cvx} = J_\mathrm{lft}\circ \Pi^{-1}$ yields
\begin{align}\label{eq:weak_PL-lft}
 J_\mathrm{lft}(K,\gamma,P)-J^\star
 \leq
\|V\|_F
 \times
 \beta_{(K,\gamma,P)},
\end{align}
where
$V\in \partial
J_\mathrm{lft}(K,\gamma,P)$
and
\begin{align*}
\beta_{(K,\gamma,P)}
=&
\sigma_\mathrm{max}
\left( \nabla\Pi^{-1}
\left(
\Pi(K,\gamma,P)
    \right)\right)\\
    &\times
        \left\|\Pi(K,\gamma,P)
    -\Pi(K^\star,J^\star,P^\star)\right\|_F>0.
\end{align*}
Note that $\beta_{(K,\gamma,P)}$ is positive for $(K,\gamma,P)\in \mathcal{L}_\mathrm{lft}.$
Recalling the partial minimization property of $J(K)$ in \cref{eq:Jk_parmin},
we now have
the following fact from \cite[Th. 10.13]{rockafellar2009variational}:
\begin{equation*}
    \partial J(K)
     \subset 
     \!\!\!
     \bigcup_{(\gamma,P)\in  {\arg\inf}_{\gamma,P} 
     J_{\mathrm{lft}}(K,\gamma,P)
     }
     \mathcal{M}(\gamma,P),
\end{equation*}
where
\begin{equation*}
    \mathcal{M}(\gamma,P)
    =\left\{
    G\middle|
    (G,0,0)\in
    \partial J_\mathrm{lft}(K,\gamma,P)
    \right\}.
\end{equation*}
Note that the assumptions in
\cite[Th. 10.13]{rockafellar2009variational}
are satisfied
due to the boundedness of the level set;
whenever
$\|K\|_F\to\infty$ or $\|P\|_F\to \infty$,
it leads to the infeasibility, and $J_\mathrm{lft}(K,\gamma,P)\to \infty$, which implies the level set 
$\left\{(K,\gamma,P) \in \mathcal{L}_\mathrm{lft} \mid J_\mathrm{lft}(K,\gamma,P)\leq t \right\}$
for any $t\geq0$ is bounded.

Then substituting $(\gamma_K,P_K)\in 
 {\arg\inf}_{\gamma,P} J_{\mathrm{lft}}(K,\gamma
 ,P)$ into \cref{eq:weak_PL-lft}
gives
$$J(K)-J^\star
\leq
\|G\|_F \times \beta_{(K,\gamma_K,P_K)},$$
where 
$G$ can be taken as any subgradient in $\partial J(K)$
thanks to the arbitrariness
of $W\in \partial_{(\gamma,Y,X)} J_\mathrm{cvx}
    \left(
    \Pi(K,\gamma,P)
    \right)$
    in \cref{eq:convexity_J_cvx}.
Thus, 
we arrive at
the weak PL inequality in \cref{eq:weakPL} 
with
\begin{align}\label{eq:mu_K}
    \mu_K=\max_{
    \substack{
    (\gamma_K,P_K)\in \arg\min_{\gamma,P}J_{\mathrm{lft}(K,\gamma,P)}}}
    \beta_{(K,\gamma_K,P_K)}^{-1}>0,
\end{align}
where $\mu_K$ is positive
as $K\in \mathcal{K}$.
\hfill$\square$
\end{document}